\documentclass[12pt]{article}
\usepackage[cp1251]{inputenc}
\usepackage[russian]{babel}
\usepackage{graphicx}
\usepackage[russian]{babel}

\newtheorem{thm}{Теорема}
\newtheorem{lemma}{Лемма}

\begin{document}

\noindent {\bf Исследование двухстаночной задачи построения расписания при наличии частичного порядка между операциями}\\

\noindent
{\sf А.А. Романова}\\

\noindent {\em Омская юридическая академия, \\Омский государственный
университет им. Ф.М. Достоевского, \\ Омск\\
e-mail: anna.a.r@bk.ru}\\

\noindent {\small {\bf Аннотация.} В работе рассматривается
NP-трудная задача составления расписания выполнения операций
единичной длительности на двух станках при наличии частичного
порядка между операциями c критерием минимизации момента завершения
всех операций. Предложен приближенный алгоритм решения задачи;
построена достижимая оценка погрешности получаемого решения.
Доказана полиномиальная разрешимость задачи в случае, когда каждая
операция, выполняющаяся на первом станке, связана отношениями
предшествования ровно с двумя операциями, выполняющимися на втором
станке; разработан соответствующий алгоритм.}

\noindent {\small {\bf Ключевые слова: } кросс-докинг, расписание,
частичный порядок, приближенный алгоритм.
}\\

\section*{Введение}

Появление в последние десятилетия новых систем доставки продукции
привело к необходимости решения проблемы их оптимального
функционирования. Рассмотрим одну из таких систем, получившую в
литературе название кросс-докинг \cite{2apte}. В системе
кросс-докинг процесс приемки и отгрузки товаров происходит через
склад напрямую, без размещения в зоне долговременного хранения, в
результате чего продукция доставляется в кратчайшие сроки.
Преимуществами такого сквозного складирования являются: быстрая
доставка продукции заказчикам, сокращение складских площадей и
снижение затрат на оплату аренды складов и труд персонала.

При двухэтапном кросс-докинге партии однотипных товаров сначала
разгружаются, подвергаются распаковке и переоформлению, затем из
этих товаров формируются заказы конечных потребителей, включающие в
себя разные типы товаров.

Пусть $n$ -- число видов продукции, доставляемых $n$ транспортными
средствами. Считаем, что разгрузкой и распаковкой занимается одна
группа рабочих, поэтому одновременно не может разгружаться более
одного транспортного средства. С другой стороны, имеется $m$ заказов
от потребителей. Очевидно, что заказ не может быть сформирован, пока
не будут готовы все необходимые для него товары. Считаем, что
оформлением заказа также занимается одна бригада рабочих.

Сформулируем данную задачу в терминах теории расписаний в системе
flow-shop с двумя станками, где операция на первой машине --
разгрузка и распаковка привезенной продукции, а операция на второй
машине -- формирование заказа. Заметим, что работа над очередным
заказом не может быть начата, пока все виды продукции из этого
заказа не будут разгружены и распакованы. Другими словами, операция
на второй машине не может начаться, пока не выполнен некоторый набор
операций на первой машине. Необходимо составить расписание
обслуживания транспортных средств, доставляющих продукцию на склад,
и операций по формированию заказов таким образом, чтобы  время
завершения формирования всех заказов было минимально.

В данной работе исследуется задача с одинаковыми длительностями
выполнения операций. Такая задача возникает на практике, когда
продукция доставляется однотипными партиями, и на оформление заказов
затрачивается одинаковое время. Без ограничения общности, будем
далее считать, что длительности всех операций равны единице.

Эффективной работе в системе кросс-докинг посвящено достаточно много
работ. Лишь в немногих из них для описания рассматриваемых задач
используются термины теории расписаний. Работа \cite{3boysen} одна
из таких. В ней представлен обзор результатов по задачам в системе
кросс-докинг, предлагается их классификация.

Рассматриваемая задача $F2|CD|C_{\max}$ является NP-трудной в
сильном смысле как для произвольных, так и для единичных
длительностей операций \cite{4chenlee}. Также в \cite{4chenlee} для
решения задачи предложен алгоритм ветвей и границ. Следует отметить,
что известная полиномиально разрешимая задача Джонсона
$F2||C_{\max}$ \cite{6jonson} является частным случаем этой задачи.
Решающее правило для задачи Джонсона используется в некоторых
работах для построения приближенного решения
\cite{4chenlee,5chensong}.

\section{Постановка задачи}

Дадим  формальную постановку рассматриваемой задачи. Имеются
множества операций $V_1 = \{A_1, \ldots, A_n\}$, $V_2 =
\{B_1,\ldots, B_m\}$ и производственная линия, состоящая из двух
машин $M_1$ и $M_2$. Операции $A_i$ выполняются на машине $M_1$,
операции $B_j$ -- на машине $M_2$. Длительности операций $A_i$ и
$B_j$ равны единице. Прерывания и одновременное выполнение двух и
более операций на одной машине не допускаются. Частичный порядок
выполнения операций задан двудольным графом $G = (V, E)$, где $V =
V_1\cup V_2$ -- множество вершин, $E = \{( A_i, B_j) \ | \ A_i
\rightarrow B_j, A_i \in V_1, B_j \in V_2\}$ -- множество дуг,
задающих отношения предшествования между операциями. Здесь запись
$A_i \rightarrow B_j$ означает, что операция $B_j$ не может начаться
до завершения операции $A_i$.

Расписание выполнения операций на машинах определяется заданием
времени начала каждой операции. Пусть $s_i^A$ -- время начала
операции $A_i$, а $s_j^B$ -- время начала операции $B_j$. Расписание
является допустимым, если к началу выполнения каждой из операций
$B_j$ завершены все операции на первой машине, связанные с ней
отношением предшествования.

В рассматриваемой задаче требуется найти допустимое расписание, при
котором время завершения всех операций
$C_{\max}=\max\{\max\limits_{i=1,\ldots,n}
(s_i^A+1),\max\limits_{j=1,\ldots,m} (s_j^B+1)\}$ минимально.

Обозначим через $S_i$ множество операций, связанных отношениями
предшествования с операцией $A_i$, а через $T_j$ -- множество
операций, связанных отношениями предшествования с операцией $B_j$.
Пусть  $d_i^A=|S_i|$, $d_j^B= |T_j|$, $i = 1,\ldots, n$, $j =
1,\ldots, m.$

Так как операции $A_i$ не имеют предшествующих, их можно выполнять
на первой машине без простоев. Поэтому расписание на первой машине
можно задать последовательностью выполнения операций. При заданной
последовательности выполнения операций на первой машине расписание
выполнения операций на второй машине, при котором время завершения
всех операций минимально, можно найти за полиномиальное
время (см. п. 2). 

В п. 2 предлагается приближенный алгоритм решения с достижимой
оценкой точности. В п. 3 выделен полиномиально разрешимый случай
задачи; приведен соответствующий алгоритм.

\section{Приближенный алгоритм}

Опишем предлагаемый алгоритм $Gr$ нахождения приближенного решения
задачи $F2|p_j=1,CD|C_{\max}$. Идея алгоритма состоит в построении
последовательности выполнения операций на первой машине ``жадным''
образом. Чем больше полустепень исхода вершины первой доли графа,
тем раньше имеет смысл выполнить соответствующую операцию.

Упорядочим операции $A_i$ по невозрастанию величин $d_i^A$. При
совпадении $d_i^A$ упорядочим соответствующие операции по
невозрастанию величин $\frac{d_i^A}{\sum\limits_{j \in S_i} d_j^B}$.
Получим последовательность $\pi$. Положим $s_{\pi_i}^A=i-1$,
$i=1,\ldots,n$.

По последовательности $\pi$ определяем расписание выполнения
операций на второй машине, минимизирующее время завершения всех
операций. Вычислим $r_j=\max\{d_j^B, \max\limits_{i\in
T_j}(s_i^A+1)\}$ -- время готовности операции $B_j$ к выполнению.
Действительно, операция $B_j$ не может начаться раньше завершения
всех предшествующих ей операций на первой машине. Далее, сортируем
последовательность операций на второй машине по неубыванию $r_j$ и в
данном порядке выполняем их, соблюдая условие $s_j^B \ge r_j$.
Данный алгоритм является полиномиальным с трудоёмкостью $O(nm)$
операций.

Обозначим через $q$ наименьший номер, при котором выполняется
неравенство:

\begin{equation}
    \sum_{i=1}^q d_{\pi_i}^A > \sum_{i=1}^n d_{\pi_i}^A - m.
\end{equation}

В следующей теореме получена оценка погрешности предложенного
алгоритма.
\begin{thm} Пусть $C_{\max}^{*}$ -- длина оптимального расписания, $C_{\max}^{Gr}$  --
длина расписания, полученного алгоритмом $Gr$. Тогда справедливо
следующее неравенство:
$$\frac{C_{\max}^{Gr}}{C_{\max}^{*}} \le \frac{\max\{q+m,n\}}{\max\{m+d_{\min}^A,n+d_{\min}^B\}},$$
где $d_{\min}^A=\min\limits_{i=1,\ldots,n} d_i^A$,
$d_{\min}^B=\min\limits_{j=1,\ldots,m} d_j^B$.

\end{thm}

\noindent {\bf Доказательство:} Очевидно, что
$\max\{m+d^A_{\min},n+d^B_{\min}\}$ -- нижняя граница длины
расписания. Построим верхнюю границу длины расписания, которое может
быть получено алгоритмом $Gr$. Рассмотрим вектор
$d^B=(d_1^B,\ldots,d_m^B)$. Заметим, что операция $A_i$ вносит
вклад, равный единице, в значения $d_j^B$ для всех $j \in S_i$.
После выполнения каждой операции на первой машине пересчитываем
вектор $d^B$, вычитая по единице из соответствующих $d_j^B$. При
этом операция $B_j$ может поступать в обработку только в случае
$d_j^B=0$. Построим допустимое расписание, соответствующее
построенной алгоритмом последовательности $\pi$ выполнения операций
на первой машине. Выполним  операции $A_{\pi_{1}}$, $A_{\pi_{2}}$,
..., $A_{\pi_{q}}$ на первой машине, не выполняя операций на второй
машине. После этого в силу неравенства (1) сумма координат вектора
$d^B$ станет меньше числа $m$ операций на второй машине. Это значит,
что после выполнения $q$ операций на первой машине по крайней мере
одна операция $B_z$ на второй машине доступна к выполнению. Докажем,
что, начиная с периода времени $q+1$, операции на второй машине
можно выполнять без простоев. В соответствии с алгоритмом операции
на первой машине выполняются по невозрастанию полустепеней исхода
соответствующих вершин графа. Пусть $v$ -- наибольшая полустепень
исхода оставшихся вершин первой доли графа. Очевидно, $v\ge 1$. Во
время выполнения соответствующей операции $A_{\pi_{q+1}}$  можно
выполнить операцию $B_z$ на второй машине. Сумма координат вектора
$d^B$ после этого уменьшится на $v$, то есть станет меньше, чем
$m-v$. При этом на второй машине останется $m-1$ операций для
выполнения. Это значит, что опять по крайней мере одна операция на
второй машине доступна для выполнения во время выполнения $(q+2)$-ой
операции на первой машине. Такая же ситуация повторяется и далее.
После очередной операции на первой машине по крайней мере одна
операция на второй машине становится доступна для выполнения. Это
продолжается до тех пор, пока вектор $d^B$ не станет нулевым. Тогда
оставшиеся операции на обеих машинам могут выполниться без простоев.
Длина соответствующего расписания будет равна  $\max\{q+m, n\}$.
Таким образом, длина расписания, полученного алгоритмом $Gr$, не
превышает $\max\{q+m, n\}$. $\bigtriangleup$

Приведём класс примеров, для которых построенная оценка погрешности
достижима. Пусть $n=k+l+s$, $m=2k+s$, $k\ge l$, $s\ge 3$. Граф
отношений предшествования представлен на рисунке 1 (каждая вершина
из блока из $l$ вершин первой доли связана с каждой вершиной из
блока из $s$ вершин второй доли).

  \includegraphics[width=0.6\textwidth]{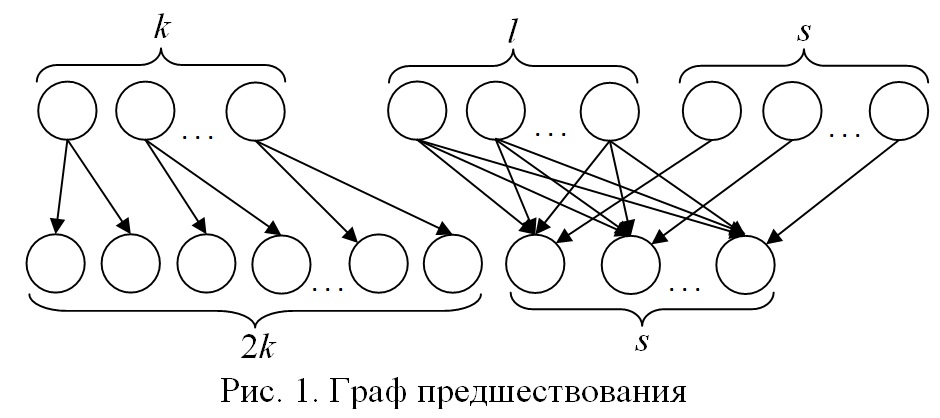}

Длина оптимального расписания равна $C_{\max}^*=2k+s+1$, что
совпадает с нижней оценкой длины расписания
$\max\{m+d^A_{\min},n+d^B_{\min}\}$. Для получения оптимального
расписания необходимо выполнить на первой машине блок из $k$
операций, затем -- блок из $l$ операций и в конце -- блок из $s$
операций. Однако в  соответствии с алгоритмом на первой машине
должен сначала выполниться блок из $l$ операций, затем -- из $k$
операций, и в конце -- из $s$ операций. В этом случае $q=l+1$. Длина
приближенного расписания равна $2k+s+l+1$, что совпадает с верхней
оценкой длины расписания $max\{q+m,n\}$.

\section{Полиномиально разрешимый случай}

Рассмотрим класс примеров  задачи $F2|p_j=1, CD|C_{\max}$, в которых
степени вершин первой доли $V_1$ графа $G$ имеют степень 2. В
дальнейшем данный класс будем обозначать через $D_2$. В [1] доказана
полиномиальная разрешимость класса $D_2$ в случае связного основания
графа $G$, построен соответствующий алгоритм $PD_2$. В данной работе
показано, что алгоритм $PD_2$ находит оптимальное решение для всех
примеров из класса $D_2$. Приведем этот алгоритм.

{\it Алгоритм  $PD_2$.}

{\it Пока} $V_2 \not= \emptyset$ выполнять:

{\it Шаг 1.} Упорядочить текущее множество $V_2$ операций  на второй
машине по неубыванию $d_j^B$. После упорядочивания имеем $d_{l_1}^B
\le d_{l_2}^B \le \ldots \le d_{l_k}^B$, где $k = |V_2|$. Операции
$B_{l_j}$, $l=1,\ldots,r$, для которых $T_{l_j}  = \emptyset$,
выполним в произвольном порядке на машине $M_2$ в момент ее
освобождения, и соответствующие вершины удалим из множества $V_2$.
Переходим на шаг 2.

{\it Шаг 2.} В момент освобождения первой машины выполним
$d_{l_{r+1}}^B$ операций  множества $T_{l_{r+1}}$. Удалим все дуги,
выходящие из этих вершин, вместе с вершинами. Выполним операцию
$B_{l_{r+1}}$ на второй машине с соблюдением условий допустимости и
удаляем соответствующую вершину из множества $V_2$. Обновляем
множества $T_j$ и величины $d_j^B$ для оставшихся во множестве $V_2$
вершин $B_j$. \\{\it Конец цикла.}

Описанный алгоритм является полиномиальным с трудоемкостью
\linebreak $O(m^2\log_2 m+mn)$ операций.

Разобьем расписание, полученное в результате  работы алгоритма, на
блоки. Пусть $\min\limits_{j=1,\ldots,m}d_j^B=t$.  Значит, первая
выбранная вершина будет степени $t$. Эта вершина вместе с вершинами,
связанными с ней отношениями предшествования, дают начало блоку
$BL_t$. Последующие включаемые в расписание операции будут отнесены
к этому же блоку до тех пор, пока на шаге 2 величина $d_{l_{r+1}}^B$
не станет больше $t$. Если в какой-то момент времени на шаге 1
станет $\min\limits_{j=1,\ldots,m}d_j^B=g>t$,  то включаемые в
расписание на шаге 2 операции дадут начало блоку $BL_g$.

Рассмотрим блок $BL_i$. Множество операций  блока $BL_i$ на первой
машине, выполняющихся до всех операций второй машины этого же блока,
будем называть \textit{отступом} блока $BL_i$, а время их выполнения
(численно равное их количеству) –- длиной отступа блока $BL_i$.
Очевидно, длина отступа в каждом таком блоке равна $i$. Множество
операций блока $BL_i$, выполняющихся на второй машине, которые не
могут начаться раньше, чем завершится последняя операция этого блока
на первой машине, будем называть \emph{выступом} блока $BL_i$, время
их выполнения –- длиной выступа. В силу структуры решаемого класса
задач длина выступа всегда равна 2. В общем виде блок $BL_i$  можно
изобразить следующим образом:

\includegraphics[width=0.6\textwidth]{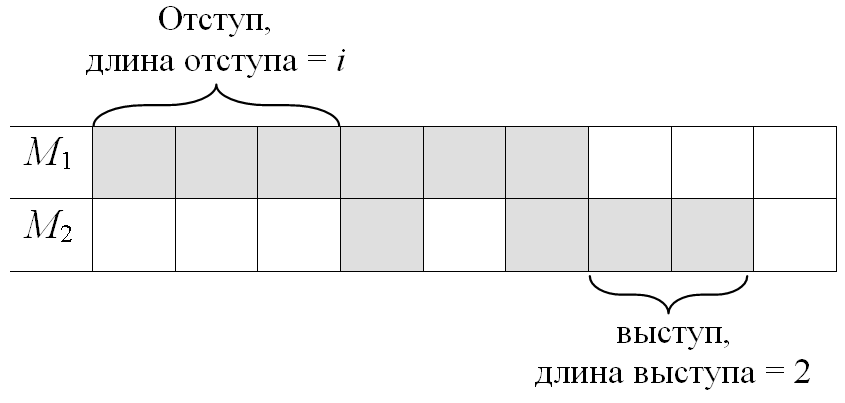}

\hspace*{2cm} {\normalsize Рис. 2. Структура блока $BL_i$.}\\[2mm]

 По построению блока $BL_i$, если на второй машине есть простои, то
продолжительность каждого из них не превышает $i$. Заметим, что блок
$BL_0$ состоит только из операций, выполняющихся на второй машине.

Расписание, получаемое алгоритмом, можно представить в виде
состыкованных блоков $BL_{v_j}$, $j=1,\ldots,h$,  при этом $v_1 <
v_2 < \ldots < v_h$. Стыковка визуально -- ``склеивание'' блоков по
операциям на первой машине так, чтобы они выполнялись без простоев.
Блок может поменять конфигурацию при стыковке. Операции на первой
машине могут сдвинуться влево из-за требования выполнения их без
простоев. Операции на второй машине могут сдвинуться вправо из-за ее
занятости в нужный момент (особенно при наличии блока $BL_0$).

 Для
иллюстрации работы алгоритма и введенных понятий и обозначений
приведем пример. Граф $G$ представлен на рис. 3.

\includegraphics[width=0.6\textwidth]{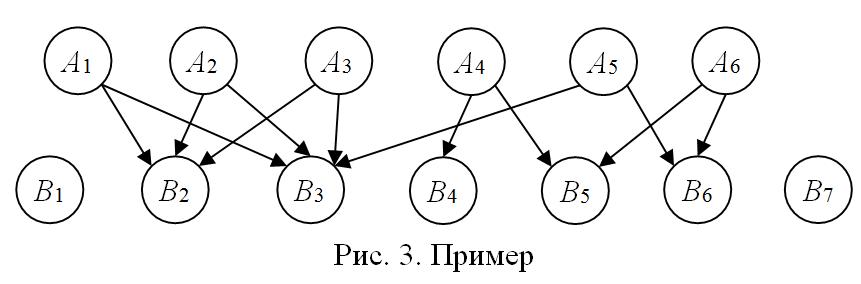}

На шаге 1 сортируем вершины второй доли графа по неубыванию числа
входящих дуг. Получим последовательность ($B_1$, $B_7$, $B_4$,
$B_5$, $B_6$, $B_2$, $B_3$) и соответствующий вектор полустепеней
захода $d^B = (0, 0, 1, 2, 2, 3, 4)$. Включаем в расписание операции
$B_1$ и $B_7$, не связанные отношениями предшествования с операциями
на первой машине, и удаляем соответствующие вершины из графа. Эти
операции образуют блок $BL_0$. На шаге 2 операции $A_4$, $B_4$ дают
начало блоку $BL_1$. Удаляем соответствующие вершины из графа вместе
с дугами, выходящими из вершины $A_4$. Переходим на шаг 1. Сортируем
оставшиеся вершины второй доли по неубыванию числа входящих дуг.
Получим $(B_5, B_6, B_2, B_3)$ и вектор $d^B = (1, 2, 3, 4)$. Так
как минимальная степень равна 1, то блок $BL_1$ продолжается. На
шаге 2 в расписание включаются операции $A_6$, $B_5$ и удаляются
дуги, выходящие из $A_6$. Переходя на шаг 1, заново сортируем
оставшиеся вершины второй доли $(B_6, B_2, B_3)$, вектор $d^B = (1,
3, 4)$. Блок $BL_1$ продолжается операциями $A_5$ и $B_6$. Остаются
вершины $B_2$, $B_3$ со степенями, равными 3. Продолжая по аналогии,
получаем, что оставшиеся вершины образуют блок $BL_3$. Блоки и
расписание, полученное в результате стыковки блоков, представлены на
рис. 4 (а, б).

\includegraphics[width=0.8\textwidth]{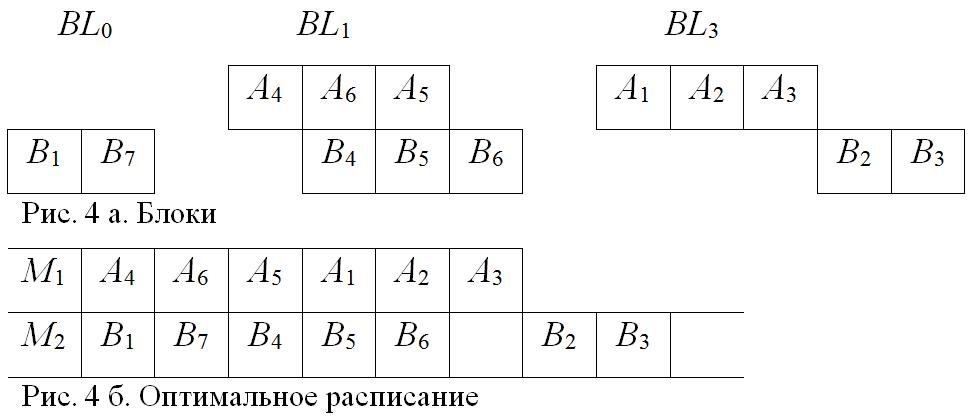}

Докажем, что описанный алгоритм находит оптимальное решение в
рассматриваемом классе задач. Начнем с двух вспомогательных лемм.

\begin{lemma}
Нижняя граница длины расписания для задач класса $D_2$ равна
$\max\{n+2,m\}$, если во второй доле графа предшествования есть
висячие вершины, и $\max\{n+2,m+1\}$, если таких вершин нет.
\end{lemma}
{\bf Доказательство:} Длина расписания в обоих случаях не меньше,
чем $n + 2$, так как каждая операция на первой машине связана ровно
с двумя операциями предшествования на второй машине, и после
завершения последней операции на первой машине должны выполниться
еще как минимум две на второй машине. Также в обоих случаях длина
расписания не может быть меньше нагрузки второй машины, то есть
величины $m$. Если во второй доле графа нет висячих вершин, то
границу можно уточнить. Полустепень захода вершин второй доли не
меньше 1, поэтому в первый момент времени машина $M_2$ будет
простаивать, и нагрузка этой машины не может быть меньше $m + 1$.
$\bigtriangleup$

\begin{lemma} В блоке $BL_k$ при $k \ge 2$ длина выступа равна либо 1,
либо 2.
\end{lemma}
{\bf Доказательство:} Рассмотрим процесс образования блока $BL_k$,
где $k \ge 2$, начиная с рассмотрения некоторой вершины  $B_{i_1}$,
полустепень захода которой равна $k$. Отметим, что в начале
формирования блока полустепени захода всех остальных вершин второй
доли не меньше $k$. Удаляя $k$ вершин первой доли графа, связанных с
вершиной $B_{i_1}$, вместе с инцидентными им $2k$ дугами, мы
уменьшаем на $2k$ сумму степеней вершин второй доли. При этом на $k$
уменьшили степень вершины $B_{i_1}$, значит, на остальные вершины
приходится суммарное уменьшение, равное $k$. После удаления дуг и
вершины  $B_{i_1}$ не может образоваться более двух вершин,
полустепень захода которых равна 1, так как для появления таких
вершин необходимо удалить не менее $k -1$ инцидентных им дуг. По той
же причине не может образоваться более одной висячей вершины. После
применения шага 2 текущая конфигурация блока изображена на рис. 5а.

\includegraphics[width=0.8\textwidth]{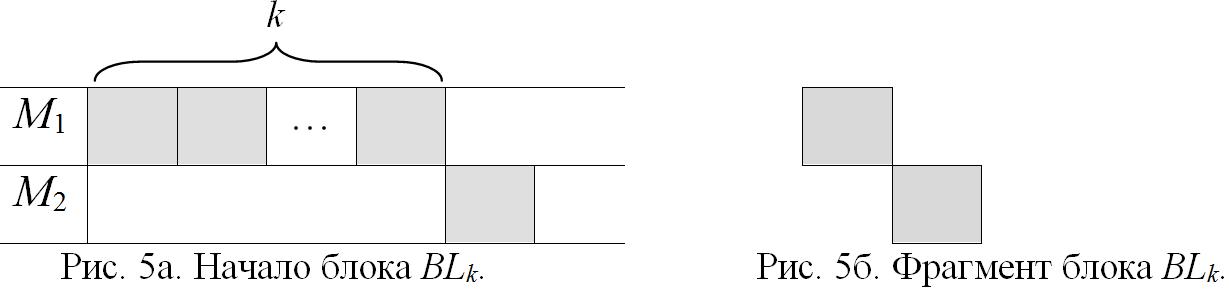}

При этом возможны следующие варианты для оставшихся вершин второй
доли:

А) минимальная степень вершин равна $t$, где $t\in \{2,\ldots, k\}$;

Б) две вершины -- степени 1, остальные -- большей степени (при их
наличии);

В) одна вершина -- степени 1, остальные -- большей степени (при их
наличии);

Г) одна висячая вершина, остальные -- большей степени (при их
наличии);

Д) пустой граф;

Е) минимальная степень вершин больше $k$.

Заметим, что при возникновении варианта Е) процесс построения блока
завершен; текущее состояние блока удовлетворяет условиям леммы.

Если появился вариант А), то на следующей итерации к блоку добавится
фрагмент, аналогичный изображенному на рис. 5а (выступ не меньше 2 и
не превышает $k$). В дальнейшем опять возможны варианты А)-Е).

При появлении варианта Б) к блоку добавляется фрагмент, изображенный
на рис. 5б. Если при этом эти две вершины были инцидентны разным
вершинам первой доли, то после этого приходим к варианту В); если
они были связаны с одной и той же вершиной, то к варианту Г).

Если появился вариант В), то к блоку также добавляется фрагмент,
изображенный на рис. 5б. В дальнейшем, очевидно, опять приходим к
варианту А), Д) или Е).

При появлении варианта Г) добавляется соответствующая операция на
вторую машину; выступ блока становится равным 2. В дальнейшем
возможны варианты А), Д) или Е).

Заметим, что при появлении варианта А) отступ очередного фрагмента
не меньше 2, поэтому при простой стыковке фрагментов не образуется
простоя на первой машине, а это значит, что выступ опять станет
равным 1. Всевозможные переходы между вариантами изображены на рис.
6.

\includegraphics[width=0.4\textwidth]{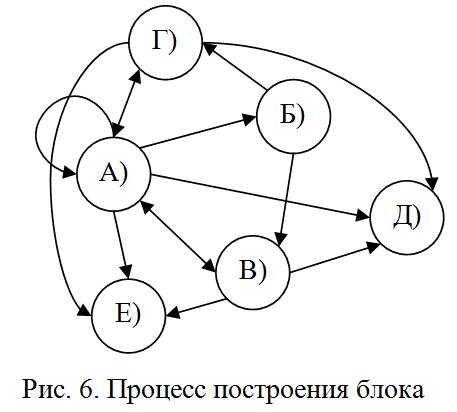}

Очевидно, процесс построения блока конечен и завершится вариантом Д)
или Е). Фрагменты, добавляемые к блоку при возникновении вариантов
А)-Г), не образуют простоя на первой машине, так как отступ любого
добавляемого фрагмента не меньше двух, а длина выступа, как видно,
всегда получается не больше двух. Так как длина отступа в фрагментах
всегда не меньше 2, а длина выступа не превышает двух, то стыковкой
получим блок с длиной выступа, не превышающей 2. $\bigtriangleup$

 Переходим к доказательству
основного результата параграфа.

\begin{thm} В задаче $F2| p_j = 1; CD|C_{\max}$, в которой каждая
операция, выполняющаяся на первой машине, связана отношениями
предшествования ровно с двумя операциями, выполняемыми на второй
машине, алгоритм $PD_2$ находит оптимальное расписание.
\end{thm}
{\bf Доказательство:} В силу леммы 2 у блоков $BL_t$, $t \ge 2$,
выступ равен 1 или 2. Будем называть такие блоки \emph{старшими}. У
последнего блока выступ обязательно равен 2 (так как последняя
операция на первой машине связана равно с двумя операциями на второй
машине). Если общее число операций на первой машине в старших блоках
равно $w$, то длина фрагмента, образованного их стыковкой, равна $w
+ 2$. Действительно, длина отступа каждого следующего блока больше,
чем длина выступа предыдущего, поэтому операции на первой машине
выполняются без простоев; выступ последнего блока равен 2.

Если в расписании отсутствуют блоки $BL_0$ и $BL_1$, алгоритм найдет
расписание длины $n + 2$. Этот случай, очевидно, соответствует
случаю отсутствия висячих вершин во второй доле графа и выполнению
неравенства $n + 2 > m$. В силу леммы 1 нижняя граница длины
расписания составит $\max\{n+2,m+1\}=n+2$.  Таким образом, найденное
расписание оптимально.

Пусть в расписании имеется блок $BL_1$, но нет висячих вершин во
второй доле (то есть блок $BL_0$ отсутствует). Длина отступа блока
равна 1. Рассмотрим процесс построения блока по аналогии со старшими
блоками. Блок начинается с некоторой операции $B_v$ с фрагмента,
изображенного на рис. 5б, после включения в блок операции на первой
машине, связанной с $B_v$, и самой операции $B_v$. После удаления
двух дуг в соответствии  с алгоритмом (одной -- к вершине $B_v$)
может образоваться либо А) одна висячая вершина во второй доле
(тогда длина выступа станет равной 2), либо Б) вершина полустепени
захода 1 (тогда процесс образования блока продолжается, и стыковка
следующего фрагмента с текущим произойдет без простоев на каждой из
машин), либо В) полустепень захода остальных вершин станет больше 1
(и тогда процесс построения блока закончен). После варианта А) мы
переходим либо в варианту Б), либо к варианту В). При этом длина
выступа увеличивается на 1. При варианте Б) длина выступа не
меняется, так как добавляется 1 операция на машину $M_1$ и одна --
на машину $M_2$. В любом случае, можно добиться выполнения операций
блока на машинах без простоев между операциями (для машины $M_2$ --
без учета простоя в первый момент времени), то есть блок $BL_1$
имеет вид (рис. 7):

\includegraphics[width=0.6\textwidth]{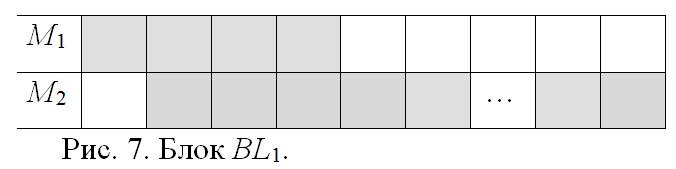}

Соединяем данный блок со старшими блоками. Выступ первого блока
может быть больше отступа первого старшего блока. В этом случае
операции старших на первой машине сдвигаем влево так, чтобы они
выполнялись без простоев. С помощью выступа первого блока сдвигаем
вправо некоторые операции второй машины, избавляясь от некоторых
простоев в старших блоках. Если $m + 1 \ge  n + 2$, то от всех
простоев на второй машине, кроме первого момента, удастся
избавиться, и длина расписания будет равна $m + 1$. Если $n + 2
> m+1$, то $(n + 1 -  m)$ простоев все же останется. И результирующее
расписание, начиная с некоторой операции, будет совпадать с
фрагментом, полученным стыковкой старших блоков, поэтому длина
выступа будет такой же, как в последнем старшем блоке, то есть 2.
Таким образом, длина расписания, найденного алгоритмом, будет равна
$n + 2$. В обоих случаях длина расписания совпадает с нижней
границей длины расписания при отсутствии висячих вершин, а, значит,
алгоритм строит оптимальное расписание.

Наконец, если в расписании присутствует блок $BL_0$ (есть висячие
вершины во второй доле графа), то при стыковке блоков $BL_0$ и
$BL_1$ образуется фрагмент (рис. 8):

\includegraphics[width=0.6\textwidth]{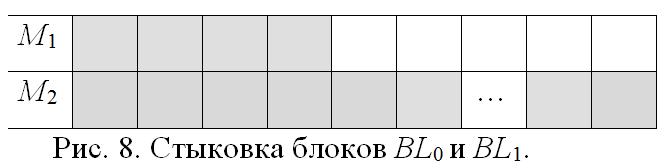}

Дальнейшие рассуждения аналогичны уже приведенным с отличием в том,
что на второй машине в первый момент нет простоя.

Если из младших блоков присутствует только блок $BL_0$, то операции
этого блока сдвигают операции старших блоков, выполняющиеся на
второй машине. Далее рассуждения проводятся по аналогии. Алгоритм
построит расписание длины $\max\{n+2,m\}$, что совпадает с нижней
границей при наличии висячих вершин. $\bigtriangleup$\\

Исследование выполнено при финансовой поддержке РФФИ (проект 12-01-00122).

\end{document}